\documentclass[12pt,a4paper]{article}
\usepackage{theorem}
\usepackage{amssymb}

\newtheorem{theorem}{Theorem}[section]
\newtheorem{proposition}[theorem]{Proposition}
\newtheorem{lemma}[theorem]{Lemma}
\newtheorem{corollary}[theorem]{Corollary}
\newtheorem{definition}[theorem]{Definition}

{\theorembodyfont{\rmfamily}\newtheorem{example}[theorem]{Example}}
\newcommand{\N}{\mathbb{N}}
\newcommand{\Z}{\mathbb{Z}}
\newcommand{\Q}{\mathbb{Q}}
\newcommand{\C}{\mathbb{C}}
\newcommand{\F}{\mathbb{F}}
\newcommand{\Mat}{{\rm Mat}}

\usepackage{hyperref}

\begin{document}

\title{Multimagic Squares}
\author{Harm Derksen,
Christian Eggermont,
Arno van den Essen}
%\date{\today}

\maketitle
\begin{abstract}
\noindent In this paper we give the first method for constructing $n$-multimagic squares (and hypercubes) for any $n$.
We give an explicit formula in the case of squares and an effective existence proof in the higher dimensional case.
Finally we prove that $n$-multimagic squares do not exist for certain orders.
\end{abstract}

\section*{Introduction}
Magic squares have been studied over $4000$ years. Recently some exciting new results have been found considering these squares. 
For instance the first method of constructing all most-perfect magic squares and their enumeration appeared in \cite{bree} (See also \cite{stewart}).
Also the non-existence of a $8\times8$ magic knights tour (\cite{magictours}) has been verified by computers.
Finally a very natural connection has been found between the properties and construction of Franklin's squares and his magic circles (\cite{pasles}).
However, there are still many unsolved problems (e.g. \cite{Abe}, \cite{Pickover}) and
as Clifford Pickover says {\sl ``the field of magic square study is wide open''} (\cite{Pickover}, p. 26).
In this paper we will concentrate on so called multimagic squares.

Suppose that $M$ is an $m\times m$ matrix $M$ consisting of natural
numbers.
Then $M$ is called a magic square if the sum of all elements in each column,
row and main diagonal gives the same number; the so-called magic number.
Let $M^{*d}$ be the matrix obtained by raising
each element of $M$ to the $d$-th power.
The matrix $M$ is called an $n$-multimagic square (where $n$ is a fixed
positive natural number) if $M^{*d}$ is a magic square
for $d=1,2,\dots,n$.
The matrix $M$ is called normal if its matrix elements consist of the consecutive
integers $1,2,\ldots,m^2$. Throughout this paper we always consider normal
magic squares (of course if $m > 1$ and $d > 1$ then the matrix $M^{*d}$ is not
normal).

\begin{figure}[ht]
$$\begin{array}{|cccccccc|} \hline
56 & 34 &  8 & 57 & 18 & 47 &  9 & 31\\ 
33 & 20 & 54 & 48 &  7 & 29 & 59 & 10\\
26 & 43 & 13 & 23 & 64 & 38 &  4 & 49\\
19 &  5 & 35 & 30 & 53 & 12 & 46 & 60\\
15 & 25 & 63 &  2 & 41 & 24 & 50 & 40\\
 6 & 55 & 17 & 11 & 36 & 58 & 32 & 45\\
61 & 16 & 42 & 52 & 27 &  1 & 39 & 22\\
44 & 62 & 28 & 37 & 14 & 51 & 21 & 3 \\ \hline
\end{array}
$$
\caption{G. Pfeffermann, $1891$}
\label{pfeff1}
\end{figure}
The first 2-multimagic square was published by Pfeffermann in 1891: it has
order 8 (Figure \ref{pfeff1}, \cite{Pfeffermann}, \cite{webBoyer}).
In 1905 the first 3-multimagic square was constructed by Tarry: it has
order 128. In 2001 both a 4- and a 5-multimagic square were constructed by
Boyer and Viricel respectively of order 512 and 1024 (\cite{Boyer}, \cite{webBoyer}
where they also give a nice history on the subject).
The record up to now was a 6-multimagic square of order 4096 constructed by Pan Fengchu in October 2003 (\cite{Fengchu}).

In this paper we give a constructive procedure to make a large 
class of $n$-multimagic squares for each positive integer $n \ge 2$. The problem of finding such
squares is reduced to an easy linear algebra problem which is solved 
in general in section~\ref{section_4}. A more explicit solution is described in section~\ref{section_5}.
This solution is used to give an explicit formula for $n$-multimagic squares for all $n\ge 3$. In
particular it gives the first 7-multimagic squares, of order $13^7$ and
8-multimagic squares of order $17^8$ etc.

The method described for constructing $n$-multimagic squares can easily be
extended to $n$-multimagic cubes and hypercubes. We refer to section~\ref{section_4} for all
definitions and more details.

\section{Preliminaries}\label{section_1}
Throughout this paper $R$ denotes a finite ring with $q$ elements.
Let $R^\star$ denote the set of elements in $R$ which have a multiplicative
inverse.
\begin{definition}\label{def1.1}
 For $c \in
R$ we call a bijection $N : R \to \{ 0,1,\ldots,q-1\}$ of type $c$ if
$ N(a) + N(-a+c) = q-1$ for all $a \in R$.
\end{definition}

\begin{lemma}\label{lem_1.1}\mbox{ }
\begin{enumerate}
\renewcommand{\theenumi}{\roman{enumi})}
\item If $2\in R^\star$, then for every $c \in R$ there exists a bijection $N$
of type $c$.
\item If $2\notin R^\star$, then for every $c \in R^\star$ there exists a
bijection $N$ of type $c$.
\end{enumerate}
\end{lemma}
{\bf Proof.} 
For $c\in R$ define $\varphi = \varphi_c : R \to R$ by
$\varphi(a) = -a + c$ for all $a \in R$. Then $\varphi^2$ is equal
to the identity.
So all orbits of $\varphi$ have length $1$ or $2$. An element $a$ is a fixed point of 
$\varphi$ if
and only if $2a = c$. For $a\in R$ denote its orbit under $\varphi$ by $O(a)$.\\
i) Let 2 be a unit in $R$. Then $\varphi$ has exactly one fixed point, namely
$a_0 := 2^{-1}c$. The orbit $O(a_0)$ of $a_0$ has one element.
Let $O(a_1),O(a_2),\dots,O(a_s)$ be the other orbits.
They all have two elements. In particular, we get $s=(q-1)/2$.
Define $N (a_0) = s$, $N(a_i) = i-1$ and $N(\varphi(a_i)) =
q-1 -N(a_i)$ ($= q-i$), for $i=1,2,\dots,s$. Now $N$ is as desired.\\
ii) Let 2 not be a unit in $R$. Then for $c \in R^\star$, $\varphi = \varphi_c$ has
no fixed points. Indeed, if $\varphi (a) = a$, then we get $2a = c\in R^\star$. 
So $2\in R^\star$ and we have a
contradiction. We can write $R = \bigcup_{i=1}^s O(a_i)$,
where each $O(a_i)$ is an orbit with two elements. In particular, $s=q/2$.
Define $N(a_i) = i-1$ and $N(\varphi(a_i)) = q-1-N(a_i)$ ($= q-i$) 
for
all $i$. Then $N$ is as desired. \hfill $\square$

\medskip\noindent
Let $m$ be a positive integer.
 For each $1 \le j \le m$ we choose a bijection $$N_{(j)} : R
\to \{ 0,1,\ldots,q-1 \}$$ of type $c_j$, 
for some $c_j \in R$ 
(this is possible
by Lemma~\ref{lem_1.1}). Put $c = (c_1,\ldots,c_m) \in R^m$ and define $N_m : R^m
\to \{ 1,2,\ldots , q^m \}$ by
$$N_m (a_1,\ldots , a_m) = 1 + \sum_{j=1}^m q^{j-1} N_{(j)}
(a_j).$$
Since the coefficients of the $q$-adic expansion of any natural number are
unique and each $N_{(j)}$ is a bijection, it follows that $N_m$ is a bijection.
\begin{lemma}\label{lem_1.2}
$N_m (-a) = q^m + 1 - N_m (a+c)$, for all $a = (a_1,\ldots,a_m) \in R^m$.
\end{lemma}
{\bf Proof.} From definition~\ref{def1.1} follows that
$$
N_m(-a)+N_m(a+c)=2+\sum_{j=1}^m q^{j-1}(N_{(j)}(-a_j)+N_{(j)}(a_j+c_j))=
$$
$$
=2+(1+q+\dots+q^{m-1})(q-1)=
2+(q^m-1)=q^m+1.
$$
%$N_m (-a) = 1 + \sum\limits_{j=1}^m q^{j-1} N_{(j)} (-a_j) = 1 +
%\sum\limits_{j=1}^m q^{j-1} (q-1 - N_{(j)} (a_j + c_j)) = 1 + \sum\limits_{j=1}^m q^j -
%\sum\limits_{j=1}^m q^{j-1} - \sum\limits_{j=1}^m q^{j-1} N_{(j)} (a_j + c_j) =
%(q^m + 1)- (1 + \sum\limits_{j=1}^m q^{j-1} N_{(j)} (a_j + c_j)) = q^m + 1 - N_m (a+c)$.
\hfill $\square$

\medskip\noindent
To conclude this section we will give a result (proposition \ref{prop_1.4})
which plays a crucial role in the next section. First some notations. 
Let $n$ and $s$ be positive integers.
Suppose that $L : R^n \to R^s$ is an affine map, i.e.,
there exists an $R$-module homomorphism $L_0:R^n\to R^s$ and a vector $v\in R^s$
such that
$$
L(a)=L_0(a)+v,\quad a\in R^n.
$$
\begin{lemma}\label{lem_1.3}
If $L : R^n \to R^s$ is a surjective 
affine map, then $\# L^{-1} (y) = q^{n-s}$ for all $y \in
R^s$.
\end{lemma}
{\bf Proof.} Let $y \in R^s$. Since $L$ is surjective there exists $a_0 \in R^n$
such that $L(a_0) = y$. It follows that $L^{-1}(y) = a_0 + \ker L_0$. Since $L$ is
surjective, so is $L_0$. It follows that $R^n / \ker L_0 \simeq R^s$,
whence $\# \ker L_0 = q^{n-s}$  and consequently $\# L^{-1} (y) = q^{n-s}$.
\hfill $\square$

\begin{proposition}\label{prop_1.4}
Suppose that $L:R^n\to R^s$ is a surjective affine map.
For each $j\in \{1,2,\dots,s\}$, let $N_{(j)} : R \to \{
0,1,\ldots, q-1 \}$ be a bijection.
For any nonnegative integers
$e_1,\ldots , e_s$ we have
$$\sum_{a\in R^n} N_{(1)} (L(a)_1)^{e_1} \ldots N_{(s)} (L(a)_s)^{e_s} = q^{n-s}
\left ( \sum_{i=0}^{q-1} i^{e_1}\right ) \ldots \left ( \sum_{i=0}^{q-1}
i^{e_s}\right ).$$
\end{proposition}
{\bf Proof.} Let $y = (y_1,\ldots , y_s) \in R^s$. Then for each $a\in
L^{-1}(y)$ we get
$$ N_{(1)}(L(a)_1)^{e_1} \ldots N_{(s)}(L(a)_s)^{e_s} = N_{(1)}(y_1)^{e_1}\ldots
N_{(s)}(y_s)^{e_s}. $$
So by Lemma \ref{lem_1.3} we obtain
\begin{equation}\label{eq_1.3}
\sum_{a\in L^{-1}(y)} N_{(1)} (L(a)_1)^{e_1}\ldots N_{(s)}(L(a)_s)^{e_s} =
q^{n-s}N_{(1)}(y_1)^{e_1}\ldots N_{(s)} (y_s)^{e_s}.
\end{equation}
Since $L$ is surjective, $R^n$ is the disjoint union of the fibres
$L^{-1}(y)$, $y\in R^s$. 
We deduce from (\ref{eq_1.3}) that
\begin{eqnarray*}
\lefteqn{\displaystyle \sum_{a\in R^n} N_{(1)} (L(a)_1)^{e_1} \ldots N_{(s)} (L(a)_s)^{e_s} \hfill}  \\
 & = & {\displaystyle \sum_{y\in R^s} \sum_{a\in L^{-1}(y)} N_{(1)}(L(a)_1)^{e_1}\ldots N_{(s)} (L(a)_s)^{e_s}} \\
 & = & {\displaystyle \sum_{y\in R^s}^{} q^{n-s} N_{(1)} (y_1)^{e_1}\ldots N_{(s)} (y_s)^{e_s}}\\
 & = & {\displaystyle q^{n-s} \left ( \sum_{y_1\in R} N_{(1)} (y_1)^{e_1} \right) \ldots \left ( \sum_{y_s\in R} N_{(s)} (y_s)^{e_s}\right )}.
\end{eqnarray*}
Since each $N_{(j)} : R \to \{ 0,1,\ldots , q-1\}$ is a bijection we get that
$$\sum\limits_{y_j \in R} N_{(j)} (y_j)^{e_j} = \sum\limits_{i=0}^{q-1}
i^{e_j},$$
which concludes the proof. \hfill $\square$

\section{A construction of n-multimagic squares}\label{section_2}
Let $n \in \N$. The following theorem gives the main tool for constructing
$n$-multimagic squares. We will use the notations introduced in the previous
section. Let $R$ be a finite ring  with $q$ elements. Write $Gl_m(R)$ for the ring of $m\times m$ invertable matrices over $R$.

First we
choose $c_1,\ldots,c_n$ in $R$ and $n$ bijections 
$$N_{(1)},\ldots , N_{(n)} : R
\to \{ 0,1,\ldots , q-1 \}$$ of types $c_1,\ldots , c_n$ respectively. 
With these (not necessarily different)
bijections we define 
$$N_n : R^n \to \{ 1,2,\ldots , q^n \}$$ as described in
section~\ref{section_1}.
This choice will be fixed throughout this section. In a similar way
we fix a bijection $$N'_{2n} : R^{2n} \to \{ 1,2,\ldots , q^{2n}\}$$
using bijections $N'_{(1)},\ldots,N'_{(2n)}$ of type $c'_1,\ldots,c'_{2n}$ respectively.

\begin{definition}
A matrix $X\in Gl_{2n}(R)$ is called an $n$-multimagic generator matrix
if for the matrices $A,B \in \Mat_{2n,n}(R)$ such that $X = (A\ B)$, all $n\times n$ minors of $A,B,A+B$ and $A-B$ are units in $R$.
\end{definition}
\begin{theorem}\label{theo_2.1}
Suppose that $X \in Gl_{2n}(R)$ is an $n$-multimagic
generator matrix.
For any $t \in R^{2n}$ the $q^n\times q^n$ matrix $M$ defined by
$$M_{N_n(a),N_n(b)} = N'_{2n}\left (X {a \choose b} +t \right ),
\quad a,b\in R^n$$
is $n$-multimagic.
\end{theorem}
Note that the matrix $M$ in the theorem is well-defined
because $N_n(a)$ takes
each value in $\{1,2,\dots,q^2\}$ exactly once,
and so does $N_n(b)$.\\[12pt]
{\bf Proof.} \\
i) First observe that all matrix elements $M_{ij}$ are distinct
since $X\in Gl_{2n}(R)$ and $N'_{2n}$
is a bijection. Consequently the
matrix $M$ consists exactly of all elements of the set $\{ 1,2,\ldots ,
q^{2n}\}$, so $M$ is normal.\\
ii) Now let $d$ be an integer with $1 \le d \le n$
and write $X= (A\ B)$ with $A,B\in\Mat_{2n,n}(R)$.
First
we want to show that the sum of all elements
in any column of $M^{*d}$ is the same constant which only
depends on $q$ and $n$.
Fix $b \in R^n$. Then the $N_n (b)$-th column of
$M^{*d}$ consists of the elements $M^d_{N_n(a),N_n(b)}$, where $a$ runs through
$R^n$ (remember that $N_n : R^n \to \{ 1,2,\ldots , q^n\}$ is a 
bijection).\\
Let $S_b(d)$ be the sum of the elements of the $N_n(b)$-th column
of $M^{*d}$, so
$$S_b(d) = \sum_{a\in R^n} M^d_{N_n(a),N_n(b)}.$$
To compute $S_b(d)$ first observe that the $j$-th component of the vector
$X{a \choose b} = (A\ B){a \choose b}$ is equal to $A_{(j)}\cdot a+B_{(j)}\cdot b$, where $A_{(j)}$ (respectively
$B_{(j)}$) denotes the $j$-th row of $A$ (respectively $B$).
Using the definitions of
$M_{N_n(a),N_n(b)}$ and $N'_{2n}$
we get
\begin{equation}\label{eq_5}
S_b(d) = \sum_{a\in R^n} \left ( 1 + \sum_{j=1}^{2n} C_j(a,b) \right )^d
\end{equation}
where
\begin{equation}\label{eq_6}
C_j(a,b) = q^{j-1} N'_{(j)} (A_{(j)}.a + B_{(j)}.b + t_j), \mbox{ for all } 1 \le
j \le 2n.
\end{equation}
Now observe that $( 1 + x_1 + \cdots + x_{2n})^d$ can be written as $1+g$, where
$g$ is a sum of terms of the form $\alpha x_{j_1}^{e_1} \ldots x_{j_s}^{e_s}$,
where $1 \le j_1 < j_2 \ldots < j_s \le 2n$, $e_1 ,\ldots , e_s \ge 1$ and $e_1
+ \cdots + e_s \le d$ (so in particular $s \le d \le n)$ and $\alpha$ is  a
positive integer. So it follows from (\ref{eq_5}) that $S_b(d)$ only depends on
$q$ and $n$ if we can show that for each set of exponents $e_1,\ldots , e_s$ and
indices $j_1,\ldots , j_s$ as above, the sum
\begin{equation}\label{eq_7}
\sum_{a\in R^n} C_{j_1} (a,b)^{e_1}\ldots C_{j_s} (a,b)^{e_s}
\end{equation}
only depends on $q$ and $n$ (and of course $e_1,\ldots , e_s$, $j_1,\ldots ,
j_s$). To see this we are going to use Proposition \ref{prop_1.4}.
Therefore put
$J = (j_1,\ldots , j_s)$ and define the affine map $L : R^n \to R^s$ by the
formula
$$L(a) = A_{(J)}\cdot a + B_{(J)}\cdot b + t_J$$
where $A_{(J)}$ (respectively $B_{(J)}$) is the $s\times n$ matrix with rows
$A_{(j_1)},\ldots ,A_{(j_s)}$ (respectively $B_{(j_1)} ,\ldots , B_{(j_s)}$) and $t_{(J)}$  in
the column of length $s$ with components $t_{j_1},\ldots t_{j_s}$. Since, as
observed above, $s \le n$ and all $n\times n$ minors of $A$ are units in $R$ (by
hypothesis) it follows that $L : R^n \to R^s$ is surjective. By
(\ref{eq_6}) we have $C_{j_i} (a,b) = q^{j_i-1}N'_{(j_i)} (L(a)_i)$
for all
$1 \le i \le s$. It follows from Proposition~\ref{prop_1.4} that the expression in
(\ref{eq_7}) is equal to
\begin{equation}\label{eq_8}
q^{n-s}\cdot q^{e_1(j_1-1)+\cdots + e_s(j_s-1)} \cdot \left ( \sum_{i=0}^{q-1} i^{e_i}
\right ) \ldots \left ( \sum_{i=0}^{q-1} i^{e_s} \right )
\end{equation}
which indeed only depends on $q$ and $n$, as desired.\\
iii) Interchanging the roles of $a$ and $b$ in the argument given in ii) we get
that all rowsums of $M^{*d}$ are equal to the same constant.\\
iv) Now let us compute the sum of the (main) diagonal elements of $M^{*d}$. This
sum is equal to
$$D(d) = \sum_{a\in R^n} M^d_{N_n(a),N_n(a)}.$$
To compute $D(d)$ we just repeat the arguments given in ii) with $b = a$. It then
remains to show that the expression in (\ref{eq_7}) with $b = a$ equals the
expression given in (\ref{eq_8}),
since this results in $D(d)$ being equal to the same constant as the rows.
Therefore just observe that
$C_{j_i} (a,a) = q^{j_i-1} N'_{(j_i)} (L_1(a)_i)$ for all $1 \le i \le s$, where $L_1 : R^n \to R^s$
is the affine map defined by
$$L (a) = A_{(J)}\cdot a + B_{(J)}\cdot a + t_{(J)} = (A+B)_{(J)}\cdot a + t_{(J)}$$
(recall that $J = (j_1,\ldots , j_s)$ ). Since by hypothesis all $n\times n$
minors of $A+B$ are units in $R$, it follows that $L$ is surjective. Then
using Proposition \ref{prop_1.4} again we obtain that the expression in
(\ref{eq_7}) with $b = a$ is indeed equal to the expression given in
(\ref{eq_8}).\\
v) Finally we compute the sum of all elements from the ``second'' diagonal of
$M^{*d}$. This sum $D'(d)$ is equal to
$$D'(d) = \sum_{a\in R^n} M^d_{N_n(a),q^n+1-N_n(a)}.$$
Since by lemma \ref{lem_1.2} $q^n +1 - N_n(a) = N_n (-a+c)$ we get
$$D'(d) = \sum_{a\in R^n} M^d_{N_n(a),N_n(-a+c)}.$$
Then repeating the arguments in ii) with $b$ replaced by $-a+c$ leads us to define the
affine map $L': R^n \to R^s$ by
$$L'(a) = A_{(J)}\cdot a + B_{(J)} (-a+c) + t_{(J)} = (A-B)_{(J)}\cdot a + (B_{(J)}\cdot c +
t_{(J)}).$$
The map $L'$ is surjective since all $n\times n$ minors of $A-B$ are units in
$R$. So again we find that the expression in (\ref{eq_7}) with $b$ replaced by $-a+c$ is
equal to the expression in (\ref{eq_8}), resulting in $D'(d)$ being equal to the same constant,
which completes the proof of this
theorem. \hfill $\square$

\section{Finding Generator Matrices}\label{section_5}
In order to be able to construct effectively $n$-multimagic squares by the
method described in theorem \ref{theo_2.1}, we need to show how to find a ring
$R$ and a $n$-multimagic generator matrix $X\in Gl_{2n}(R)$ 
which satisfy the conditions of that theorem. 

Below we describe an explicit construction of ($n$-multimagic) generator matrices.
For an effective existence proof in a more general setting see section~\ref{section_3}.

\begin{lemma}\label{lem_5.1}
Let $n \ge 1$ and $R$ a ring such that 2 and 3 are units in $R$. 
If $A$ is an $2n\times n$ matrix such that every $n\times n$ minor
of $A$ is a unit in $R$, then there exists an $2n\times n$ matrix
$B$ such that $(A\ B)$ is an $n$-multimagic generator matrix.
More precisely, if $A={P\choose Q}$ with $P,Q\in \Mat_{n,n}(R)$,
then we can take $B={\phantom{-}2P\choose -2Q}$.
\end{lemma}
{\bf Proof.} Since 2 is a unit in $R$ the hypothesis  on $P$ and $Q$ implies
that also the $n\times n$ minors of $B$ are units in $R$. Furthermore the
$n\times n$ minors of 
$A+B = {{3P} \choose {-Q}}$
are also units in $R$, since 3 and $-1$ are.
Similarly the $n\times n$ minors of $A-B$ are units in $R$. Finally, using
elementary column operators one can reduce the matrix
$(A\ B) = {{P\ \phantom{-}2P} \choose {Q\ -2Q}}$
to the matrix
${{P\ \phantom{-4}0} \choose {Q\ -4Q}}$
which is clearly invertible over $R$ since
both $\det P$ and $\det (-4Q)$ are units in $R$. \hfill $\square$

\begin{lemma}\label{lem_5.2}
Let $n \ge 2$ and define the $2n\times n$ matrix $A$
 by
$A_{i,j}=(i-1)^{j-1}$ for $i=1,2,\dots,2n-1$ and $j=1,2,\dots,n$
$(0^0=1)$,
$A_{2n,j}=0$ for $j=1,2,\dots,2n-1$ and $A_{2n,n}=1$.
So we have
$$
A=\left(\begin{array}{ccccc}
1 & 0 & 0 & \cdots & 0\\
1 & 1 & 1 & \cdots & 1\\
1 & 2 & 4 & \cdots & 2^{n-1}\\
1 & 3 & 9 & \cdots & 3^{n-1}\\
\vdots & \vdots &\vdots & & \vdots\\
1 & (2n-2) & (2n-2)^2 & \cdots & (2n-2)^{n-1}\\
0 & 0 & 0 & \cdots & 1
\end{array}\right).
$$
If $R$ is a ring then we may view $A$ as a matrix with entries in $R$.
If $\{1,2,\dots,2n-2\}\subseteq R^\star$, then all
$n\times n$ minors of $A$ are units in $R$.
\end{lemma}
{\bf Proof.} Using Vandermonde determinants one easily verifies that each 
factor
appearing in each $n\times n$ minor of $A$ is  of the form $i-j$
where $0 \le j < i \le 2n-2$, from which the desired result follows. \hfill
$\square$

\medskip\noindent
As an immediate consequence of the lemmas above we get

\begin{corollary}\label{cor_5.3}
Let $n\geq 2$ and $R$ be a ring such that $\{3,1,2,\dots,2n-2\}\subseteq
R^\star$, then the matrix $(A\ B)$ is an $n$-multimagic generator matrix
where $A$ is as in Lemma~\ref{lem_5.2} and $B$ is as in Lemma~\ref{lem_5.1}.
\end{corollary}

In theorem~\ref{cor_5.4} below we choose one bijection $N : R \to \{
0,1,\ldots , q-1 \}$ of some type $c \in R$ and define for each $m \ge 1$
$$N_m (a_1,\ldots , a_m) = 1 + \sum_{j=1}^m q^{j-1} N(a_j).$$
So in the definition of $N_m$ as given in section~\ref{section_1} we take all $N_{(j)}$ to be
equal to $N$.

\begin{theorem}\label{cor_5.4}
(An explicit formula for $n$-multimagic squares.)\\
Let $n \ge 3$, $q$ a prime number $\ge 2n-1$, $R = \F_q$ and $N : R\to \{ 0,1,\ldots , q-1
\}$ the bijection (of type $-1$) given by $N(\overline i ) = i$ for all $0 \le i
\le q-1$. Let $A$ and $B$ be as in \ref{cor_5.3}. Then for $X=(A\ B)$ and each $t \in R^{2n}$
the matrix $M$ defined in \ref{theo_2.1} is $n$-multimagic.\\
\end{theorem}

\noindent
In other words, for every positive integer $n$ there exists a (normal) $n$-multimagic
square.
In particular for $n = 7$ we get 7-multimagic squares of orders $13^7,\ 17^7,\ 19^7$ etcetera.
For $n = 8$ we get 8-multimagic squares of orders $17^8,\ 19^8,\ 23^8$ etcetera.

\section{Multimagic cubes and hypercubes}\label{section_4}
In this section we briefly indicate how the method developed in the previous
sections can be extended to construct multimagic cubes, perfect multimagic cubes and
hypercubes.

\subsection{Multimagic cubes}
Note that there is no consensus on the definition of multimagic cubes (hypercubes etc.) in the literature.
The choice given below can also be found in \cite{wolfram}, \cite{martin}.

A cube of numbers (respectively the consecutive numbers $1,2,\ldots,n^3$) is called magic
(respectively normal magic) if the sum of all elements in
each row, column and pillar is the same and is equal to the sum of all 
elements
of each of the four space diagonals. Furthermore, if $n \ge 1$ such a cube is
called $n$-multimagic if for each $1 \le d \le n$ the cube obtained by raising
each of its elements to the $d$-th power is magic.

\medskip\noindent
Completely analogues to the construction of $n$-multimagic squares in
\ref{theo_2.1} we define a $q^n\times q^n\times q^n$ cube by the formula
$$M_{N_n(a),N_n(b),N_n(c)} = N'_{3n} \left ( (A\ B\ C)  \left ( \begin{array}{l}
a\\ b\\ c\end{array} \right ) + t\right )$$
where each of the vectors $a,b$ and $c$ runs through $R^n$, $t \in R^{3n}$ and
$A,B$ and $C$ are matrices in $\Mat_{3n,n}(R)$ which satisfy the following
properties (which guarantee the matrix $M$ to be an $n$-multimagic cube):
\begin{itemize}
\item[1.] $(A\ B\ C) \in Gl_{3n}(R)$ (which guarantees that all the natural number
$1,2,\ldots,q^{3n} (= (q^n)^3)$ appear in $M$).
\item[2.] all $n\times n$ minors of the matrices $A,B$ and $C$ are units in $R$
(which guarantees that for each $1 \le d \le n$ the sum of all elements in each
column, row and pillar of $M^{*d}$ is the same, and hence equal to the magic
sum).
\item[3.] all $n\times n$ minors of the matrices $A+B+C$, $-A+B+C$, $A-B+C$ and
$A+B-C$ are units in $R$ (which guarantees that for each $1 \le d \le n$ the sum
of all elements on each of the four space diagonals of $M^{*d}$ is equal to the
magic sum).
\end{itemize}

\noindent Recall that a magic cube is called perfect if additionally the diagonals 
of each orthogonal slice have the magic sum property. Furthermore, if $n \ge 1$ such a
cube is called $n$-multimagic perfect if for each $1 \le d \le n$ the cube
obtained by raising each of its elements to the $d$-th power is perfect magic.

\medskip\noindent
To guarantee that a $n$-multimagic cube $M$ as
above is also
$n$-multimagic perfect we impose on the matrices $A,B,C$ the following
conditions
\begin{itemize}
\item[4.] all $n\times n$ minors of the matrices $A+B$, $A-B$, $A+C$, $A-C$,
$B+C$ and $B-C$ are units in $R$.
\end{itemize}

\subsection{More Generator Matrices}\label{section_3}

\noindent To find a ring $R$ and matrices $A,B$ and $C$ satisfying the
properties 1, 2, 3 and 4 one can use the method described below.
To facilitate generalizations to (even) higher dimensions we give a more general notion of generator matrix:

\begin{definition}
Let $d$ be a positive integer $\ge 2$.
We call a matrix $X \in Gl_{dn}(R)$ an $n$-multimagic $d$-generator matrix
if when we write $X = (A_1\ \ldots\ A_d)$ with $A_1,\ldots,A_d \in \Mat_{dn,n}(R)$, we have that
for $\delta_1,\ldots,\delta_d\in\{-1,0,1\}$ not all $0$, the matrix
$$\sum_{i=1}^d \delta_i A_i$$
has all $n\times n$ minors in $R^\star$.
\end{definition}

\noindent Note that the four properties of the matrix $(A\ B\ C)$ of the previous section are equivalent to saying that
it is a $n$-multimagic $3$-generator matrix.

\begin{example}
A 1-multimagic 3-generator matrix with $R = \F_q$ with $q$ prime and $q\ge 11$ :
$$\left ( \begin{array}{rrr}
1 & 2 & 4\\
1 & 2 & -4\\
1 & -2 & -4
\end{array} \right) $$
\end{example}

\begin{example}
A 2-multimagic 3-generator matrix with $R = \F_q$ with $q$ prime and $q\ge 11$ :
$$\left ( \begin{array}{rrrrrr}
1 & 0 & 2 & 0 & 4 & 0\\
1 & 1 & 2 & 2 & 4 & 4\\
1 & 2 & 2 & 4 & -4 & -8\\
1 & 3 & 2 & 6 & -4 & -12\\
1 & 4 & -2 & -8 & -4 & -16\\
0 & 1 & 0 & -2 & 0 & -4
\end{array} \right) $$
\end{example}

\begin{example}
A 2-multimagic 2-generator matrix with $R = \F_q$ with $q$ prime and $q\geq 11,\ q\neq 17$ :
$$\left ( \begin{array}{rrrrrr}
 1 & 0 & 0 & -1 & -2 & -1\\
 0 & 1 & 0 & -1 & -1 & -2\\
 0 & 0 & 1 & -2 & -1 & -1\\
 1 & 2 & 1 &  1 &  0 &  0\\
 1 & 1 & 2 &  0 &  1 &  0\\
 2 & 1 & 1 &  0 &  0 &  1
\end{array} \right) $$
\end{example}

\noindent
The following (well-known) lemma is the crucial tool 
in finding effectively $d$-generator matrices for all $d$.

\begin{lemma}\label{lem_3.1}
Let $m \ge 1$ and $Q(x_1,\ldots,x_m)$ a non-zero polynomial in the variables
$x_1,\ldots , x_m$ over $\Z$. Then one can effectively find integers
$a_1,\ldots,a_m$ such that $Q(a_1,\ldots , a_m)$ $ \ne 0$.
\end{lemma}
{\bf Proof.} We will use induction on $m$. 
The case $m = 1$ is obvious since $Q(x_1)$ has
at most $\deg Q$ zeros. Now let $m \ge 2$. Write $Q = q_d(x_1,\ldots ,
x_{n-1})x_m^{d} + \cdots + q_1(x_1,\ldots , x_{m-1}) x_m + q_0 (x_1,\ldots ,
x_{m-1})$ with $q_d \ne 0$. By the induction hypothesis there exist integers
$a_1,\ldots , a_{m-1}$ such that $q_d(a_1,\ldots , a_{m-1}) \ne 0$. So the
polynomial $q(x_m) = Q(a_1,\ldots , a_{m-1},x_m)$ in $\Z[x_m]$ is non-zero and
has $x_m$-degree $d$. Consequently there exists an integer $a_m$ such that
$$q(a_m)=Q (a_1,\ldots , a_{m-1},a_m) \ne 0,$$ 
as desired. \hfill
$\square$

\begin{proposition}%\label{prop_3.2}
Let $n \ge 1$ and $d \ge 2$. Then one can effectively find
a positive integer $q > 1$ and
matrices $A_1,\ldots,A_d \in \Mat_{dn,n}(\Z/q\Z) $ such that
$X = (A_1\ \ldots\ A_d) \in Gl_{dn}(\Z/q\Z)$ is
an $n$-multimagic $d$-generator matrix.
\end{proposition}
{\bf Proof.} To avoid complicating an easy matter we only give a proof for $d=2$.
For other $d$ the procedure is similar.\\
Let $A_u = (A_{i,j})$ and $B_u = (B_{i,j})$ be two universal
$2n\times n$ matrices, i.e., the entries $A_{i,j}$ and $B_{i,j}$ are distinct
variables. Then each $n\times n$ minor of $A_u$, $B_u$, $A_u+B_u$ and 
$A_u-B_u$ is
a non-zero polynomial in the $4n^2$-variable polynomial ring 
$\Z[A_{ij},B_{ij},
1 \le i \le 2n, 1 \le j \le n]$. Let $P$ be the product of all these minors and
let $Q$ be the product of $P$ and the polynomial $\det (A_u\ B_u)$. By lemma
\ref{lem_3.1} we can find integers $a_{i,j}$ and $b_{i,j}$ such that $Q(a_{i,j},
b_{i,j})$ is a non-zero integer. Finally let $q$ be a positive integer $> 1$ 
such that $\gcd(Q(a_{i,j},b_{i,j}),q) = 1$ for all $i,j$. 
Then one easily verifies
that $X = (a_{ij})$ and $B = (b_{ij})$ represent matrices in $\Mat_{2n,n}(\Z/q\Z)$
having the desired properties. \hfill $\square$

\subsection{Perfect Multimagic Hypercubes}
From the above it is now clear how to generalize these definitions and
constructions to higher dimensional hypercubes.

\begin{theorem}
For all integers $d\ge 2$ and $n\ge 1$ there exists a normal perfect $n$-multimagic $d$-dimensional hypercube.
\end{theorem}

\section{Orders}
We take a short look at the possible orders a $n$-multimagic square might have. First we will show that we can use $n$-multimagic squares to construct new $n$-multimagic squares of different orders.

\begin{definition}
Let $A\in\Mat_{m,m}(R)$ and $B\in\Mat_{n,n}(R)$. We define $A\star B\in\Mat_{mn,mn}(R)$ by $$A\star B_{mk+i,ml+j} = m^2 B_{k,l} + A_{i,j}$$
where $0\le i,j \le m-1$ and $0\le k,l\le n-1$.
\end{definition}

\begin{proposition}
If $A\in\Mat_{m,m}(R)$ and $B\in\Mat_{n,n}(R)$ are $p$-multimagic squares then $A\star B$ is an $p$-multimagic square.
\end{proposition}
{\bf Proof.} 
The proof is a straightforward calculation so we only explicitly show that the columnsum is
independent of the columnindex. We use the same notations as in the definition above.
Fix $1\le { e}\le p$.
We denote
$S_{ x} = \sum_i A_{i,j}^{ x}$ and $T_{ x} = \sum_i B_{i,j}^{ x}$ for $1\le{ x}\le p$.
Since $A$ and $B$ are $p$-multimagic squares these $S_x$ and $T_x$ are constants.
Fix the columnindex $\beta = ml+j$ and write $\alpha = mk+i$.
\begin{eqnarray*}
\lefteqn{ \sum_{\alpha=0}^{mn-1} (A\star B)_{\alpha,\beta}^{e} =
      \sum_{k=0}^{m-1} \sum_{i=0}^{n-1} \left ( m^2 B_{k,l} +A_{i,j} \right )^{e} } \\
& = & \sum_{k}    \sum_{i} \sum_{x=0}^{{ e}} {{ e}\choose{ x}} \left ( m^2 B_{k,l} \right )^{e-x} A_{i,j}^{x}
  =   \sum_{k}    \sum_{x} {{e}\choose{x}} m^{2(e-x)} B_{k,l}^{e-x} \sum_i A_{i,j}^{x} \\ 
& = & \sum_{k}    \sum_{x} {{ e}\choose{ x}} m^{2(e-x)} S_{x} B_{k,l}^{e-x}
  =   \sum_{{x}} {{e}\choose{x}} m^{2(e-x)} S_{x} \sum_{k} B_{k,l}^{e-x}  \\
& = & \sum_{{x}} {{e}\choose{x}} m^{2(e-x)} S_{x} T_{ e-x} 
\end{eqnarray*}
So the columnsum is independent of $\beta$, hence a constant.
\hfill$\square$\\

\noindent From the foregoing one might get the impression that $n$-multimagic squares are constructable for all (large enough) orders.
This is not the case as we will prove next.

\begin{definition}
If $p$ is a prime and $n$ is a positive integer,
then we define $v_p(n)$ as the largest nonnegative integer $e$
such that $p^e$ divides $n$. For a non-zero rational number $\frac{n}{m}$
we define $v_p(\frac{n}{m})=v_p(n)-v_p(m)$.
\end{definition}

\begin{lemma}\label{divlem} (basic properties of $v_p$)\\
Let $a,e,m,p\in \N$ with $p$ prime.
\begin{itemize}
\item[i)] $v_p(r \cdot s) = v_p(r)+v_p(s)$ for any non-zero rational numbers $r$ and $s$
\item[ii)] If $v_p(m)=e$ and $1\le a\le p^e$ then $v_p(m^2-a)=v_p(a)$.
\end{itemize}
\end{lemma}
{\bf Proof.}
Straightforward using the unique factorization of integers. \hfill$\square$\\

\begin{theorem} 
Let $M$ be a normal $n$-multimagic square of order $m$ and $p\in\N$ a prime.
If $v_p(m)=e\ge 1$ then $n\le p^{e+1}-2$.
\end{theorem}
{\bf Proof.} 
Write $M = (m_{i,j})_{i,j=1..m}$.
Note that if $f\in\C[x]$, $f(\Z)\subseteq\Z$ and $deg_x(f)\le n$ then the matrix
$M^{\circ f} = (f(m_{i,j}))_{i,j=1..m}\in\Mat_{m,m}(\Z)$ is also a (general) magic square.
Take $f(x) = {{x-1}\choose{n}} \in\Q[x]$ and note that it is a polynomial of degree $n$ with $f(\Z)\subseteq\Z$.
Since $M^{\circ f}$ is magic the sum of all $m^2$ elements of $M^{\circ f}$ is an integer which is $m$ times the magic sum, i.e. the sum of any row or column.
Since $f(\Z)\subseteq\Z$ the magic sum is an integer which implies that $m$ divides
$$\sum_{x=1}^{m^2}f(x) = \sum_{x=1}^{m^2} {{x-1}\choose{n}} = {{m^2}\choose{n+1}}.$$
This means that the (rational) number
\begin{equation}\label{eq_}
\frac{m(m^2-1)(m^2-2)\cdots(m^2-n)}{(n+1)!}
\end{equation}
is really an integer.\\
Since an $n$-multimagic square of order $m$ can not exist if an $(n-1)$-multimagic square of order $m$ does not exist,
it is enough to show that the case $n=p^{e+1}-1$ is impossible by showing that (\ref{eq_}) is not an integer.
\\
So assume $n=p^{e+1}-1$.
Note that (\ref{eq_}) is not equal to $0$ since $p^e$ divides $m$, so in particular $m^2\ge p^{2e}>p^{e+1}-1 = n$ (remember $e$ is positive).
Using lemma~\ref{divlem} we see that %we see that $v_p((m^2-i)/i) = 0$ for all $i=1,\ldots ,n$.

{ \renewcommand{\arraystretch}{2.5}
$ 
\begin{array}{rclr} 
\displaystyle
v_p \left ( \frac{m(m^2-1)(m^2-2)\cdots(m^2-n)}{(n+1)!} \right ) & = & &\\
\displaystyle
v_p \left ( \frac{m}{n+1} \cdot \frac{(m^2-1)}{1} \cdot \frac{(m^2-2)}{2} \cdots \frac{(m^2-n)}{n} \right )  & = & & \hspace*{20 pt}\mbox{{\small (\,use lemma~\ref{divlem} i)\,) }} \\
\displaystyle
v_p \left ( \frac{m}{n+1} \right ) + \sum_{i=1}^{n} {v_p \left ( \frac{m^2-i}{i} \right ) } & = &  & \mbox{{\small (\,use lemma~\ref{divlem} ii)\,) }} \\
\displaystyle
v_p \left ( \frac{m}{n+1} \right ) + \sum_{i=1}^n 0 & = & &\\
\displaystyle
v_p \left ( \frac{m}{p^{e+1}}\right ) & = & &\\
\displaystyle
 e - (e+1) & = & -1 &
\end{array}
$
}

\noindent So (\ref{eq_}) is not a integer if $n=p^{e+1}-1$.
\hfill$\square$\\

\noindent In particular this shows that there are no 3-multimagic squares of order $m\equiv 2\ (\mbox{mod}\: 4)$.

\section{More Examples}

\medskip\noindent
To conclude this paper we give some interesting multimagic squares using theorem~\ref{theo_2.1}.\\

In \ref{ex_5.5} - \ref{ex_5.7} below we choose one bijection $N : R \to \{
0,1,\ldots , q-1 \}$ of some type $c \in R$ and define for each $m \ge 1$
$$N_m (a_1,\ldots , a_m) = 1 + \sum_{j=1}^m q^{j-1} N(a_j).$$
So in the definition of $N_m$ as given in section~\ref{section_1}  we take all $N_{(j)}$
to be equal to $N$.
The notations are as in theorem~\ref{theo_2.1}, where we take $N'_{2n} = N_{2n}$.

\begin{example}\label{ex_5.5}
(A family of associative bimagic squares of order 16.) \\
Take $R = \F_2 [x] / (x^2 + x + 1)$,
{\small
$$X = \left ( \begin{array}{llll}
x&0&1&x\\
0&1&1&1\\
1&1&1&0\\
x&1&0&x
\end{array} \right ),
$$ }
$t\in R^4$ arbitrary and $N : R \to \{ 0,1,2,3 \}$ (a bijection of
type $x +1)$ given by $N(0) = 0$, $N(1) = 2$, $N(x) = 1$ and $N(x+1) = 3$.\\
Then the corresponding matrix $M$ defined in theorem \ref{theo_2.1} is
bimagic (= 2-multimagic) and associative (= the sum of any pair of matrix elements which are
symmetric with respect to the center of the square is equal to $16^2 +1$).
{\scriptsize
$$
\begin{array}{|cccc|cccc|cccc|cccc|} \hline
  41 & 252 &  74 & 155 & 125 & 176 &  30 & 207 & 129 &  84 & 226 &  51 & 213 &   8 & 182 & 103\\
  62 & 239 &  93 & 144 & 106 & 187 &   9 & 220 & 150 &  71 & 245 &  40 & 194 &  19 & 161 & 116\\
   3 & 210 & 100 & 177 &  87 & 134 &  56 & 229 & 171 & 122 & 204 &  25 & 255 &  46 & 160 &  77\\
  24 & 197 & 119 & 166 &  68 & 145 &  35 & 242 & 192 & 109 & 223 &  14 & 236 &  57 & 139 &  90\\ \hline
 240 &  61 & 143 &  94 & 188 & 105 & 219 &  10 &  72 & 149 &  39 & 246 &  20 & 193 & 115 & 162\\
 251 &  42 & 156 &  73 & 175 & 126 & 208 &  29 &  83 & 130 &  52 & 225 &   7 & 214 & 104 & 181\\
 198 &  23 & 165 & 120 & 146 &  67 & 241 &  36 & 110 & 191 &  13 & 224 &  58 & 235 &  89 & 140\\
 209 &   4 & 178 &  99 & 133 &  88 & 230 &  55 & 121 & 172 &  26 & 203 &  45 & 256 &  78 & 159\\ \hline
  98 & 179 &   1 & 212 &  54 & 231 &  85 & 136 & 202 &  27 & 169 & 124 & 158 &  79 & 253 &  48\\
 117 & 168 &  22 & 199 &  33 & 244 &  66 & 147 & 221 &  16 & 190 & 111 & 137 &  92 & 234 &  59\\
  76 & 153 &  43 & 250 &  32 & 205 & 127 & 174 & 228 &  49 & 131 &  82 & 184 & 101 & 215 &   6\\
  95 & 142 &  64 & 237 &  11 & 218 & 108 & 185 & 247 &  38 & 152 &  69 & 163 & 114 & 196 &  17\\ \hline
 167 & 118 & 200 &  21 & 243 &  34 & 148 &  65 &  15 & 222 & 112 & 189 &  91 & 138 &  60 & 233\\
 180 &  97 & 211 &   2 & 232 &  53 & 135 &  86 &  28 & 201 & 123 & 170 &  80 & 157 &  47 & 254\\
 141 &  96 & 238 &  63 & 217 &  12 & 186 & 107 &  37 & 248 &  70 & 151 & 113 & 164 &  18 & 195\\
 154 &  75 & 249 &  44 & 206 &  31 & 173 & 128 &  50 & 227 &  81 & 132 & 102 & 183 &   5 & 216\\ \hline
\end{array}
$$
}
\end{example}

\begin{example}\label{ex_5.6}
(A family of bimagic squares of odd order.)\\
Take $R = \Z/q\Z$ with $q \ge 3$, $q$ odd,
{\small
$$X = \left ( \begin{array}{llll}
0&1&1&0\\
2&0&0&1\\
1&1&2&1\\
2&1&2&2
\end{array}\right ),$$
}$t \in R^4$ arbitrary and $N : R \to \{ 0,1,\ldots , q-1 \}$ the standard bijection (of type $-1$) as in \ref{cor_5.4}.
Then the corresponding $q^2\times q^2$ matrices $M$ as defined in theorem
\ref{theo_2.1} are bimagic. In particular this gives a family of bimagic squares
of odd order.
{%\scriptsize
$$
\begin{array}{|ccc|ccc|ccc|} \hline
 1 & 35 & 60 & 23 & 48 & 79 & 18 & 40 & 65\\ 
70 & 14 & 39 & 56 &  9 & 31 & 78 & 19 & 53\\
49 & 74 & 27 & 44 & 69 & 10 & 30 & 61 &  5\\ \hline
38 & 72 & 13 & 33 & 55 &  8 & 52 & 77 & 21\\
26 & 51 & 73 & 12 & 43 & 68 &  4 & 29 & 63\\
59 &  3 & 34 & 81 & 22 & 47 & 64 & 17 & 42\\ \hline
75 & 25 & 50 & 67 & 11 & 45 & 62 &  6 & 28\\
36 & 58 &  2 & 46 & 80 & 24 & 41 & 66 & 16\\
15 & 37 & 71 &  7 & 32 & 57 & 20 & 54 & 76\\ \hline
\end{array}
$$
}
\noindent
Finally, if we choose $q = 3$ and
$t = \left(\begin{array}{l} 2\ 1\ 2\ 0\end{array}\right)^t$
we recover the associative $9\times 9$ bimagic square
constructed by R.V. Heath (see p. 212 \cite{RouseBall}) from before 1974.
\end{example}

\begin{example}\label{ex_5.7}
(An associative, pandiagonal, bimagic, .... magic square of order 25.)\\
Take $R = \F_5$, $N : R \to \{ 0,1,2,3,4 \}$ the standard bijection (of type $-1$) and 
{\small
$$X = \left (
\begin{array}{llll}
1&1&2&2\\ 
1&2&2&4\\
1&3&3&4\\
1&4&3&2
\end{array}\right ),\ t = \left ( \begin{array}{l} 0\\ 4\\ 0\\ 2
\end{array} \right ).$$}
Then the corresponding $25\times 25$ matrix is associative,
pandiagonal, bimagic and has the following properties:
\begin{itemize}
\item[i)] Each of the 25 standard $5\times 5$ submatrices is pandiagonal (even
with the same magic sum).
\item[ii)] For each pair $(i,j)$ $(1 \le i, j \le 25$) the $5\times 5$ matrix
obtained by deleting each row with row number not equivalent to $i$ mod 5 and
each column with column number not equivalent to $j$ mod 5 is pandiagonal!
\end{itemize}
{\tiny
$$
\begin{array}{@{}|c@{\hspace{3pt}}c@{\hspace{3pt}}c@{\hspace{3pt}}c@{\hspace{3pt}}c@{\hspace{3pt}}|@{\hspace{3pt}}c@{\hspace{3pt}}c@{\hspace{3pt}}c@{\hspace{3pt}}c@{\hspace{3pt}}c@{\hspace{3pt}}|@{\hspace{3pt}}c@{\hspace{3pt}}c@{\hspace{3pt}}c@{\hspace{3pt}}c@{\hspace{3pt}}c@{\hspace{3pt}}|@{\hspace{3pt}}c@{\hspace{3pt}}c@{\hspace{3pt}}c@{\hspace{3pt}}c@{\hspace{3pt}}c@{\hspace{3pt}}|@{\hspace{3pt}}c@{\hspace{3pt}}c@{\hspace{3pt}}c@{\hspace{3pt}}c@{\hspace{3pt}}c@{\hspace{3pt}}|@{}} \hline
103 & 350 & 567 & 164 & 381 & 291 & 513 & 235 & 452 & 74 & 584 & 176 & 423 & 20 & 362 & 147 & 494 & 86 & 308 & 530 & 440 & 32 & 254 & 621 & 218\\ 
 167 & 389 & 106 & 328 & 575 & 460 & 52 & 299 & 516 & 238 & 23 & 370 & 587 & 184 & 401 & 311 & 533 & 130 & 497 & 94 & 604 & 221 & 443 & 40 & 257\\
 331 & 553 & 175 & 392 & 114 & 524 & 241 & 463 & 60 & 277 & 187 & 409 & 1 & 373 & 595 & 480 & 97 & 319 & 536 & 133 & 43 & 265 & 607 & 204 & 446\\
 400 & 117 & 339 & 556 & 153 & 63 & 285 & 502 & 249 & 466 & 351 & 598 & 195 & 412 & 9 & 544 & 136 & 483 & 80 & 322 & 207 & 429 & 46 & 268 & 615\\
 564 & 156 & 378 & 125 & 342 & 227 & 474 & 66 & 288 & 510 & 420 & 12  & 359 & 576 & 198 & 83 & 305 & 547 & 144 & 486 & 271 & 618 & 215 & 432 & 29\\ \hline
 134 & 476 & 98 & 320 & 537 & 447 & 44 & 261 & 608 & 205 & 115 & 332 & 554 & 171 & 393 & 278 & 525 & 242 & 464 & 56 & 591 & 188 & 410 & 2 & 374 \\
 323 & 545 & 137 & 484 & 76 & 611 & 208 & 430 & 47 & 269 & 154 & 396 & 118 & 340 & 557 & 467 & 64 & 281 & 503 & 250 & 10 & 352 & 599 & 191 & 413\\
 487 & 84 & 301 & 548 & 145 & 30 & 272 & 619 & 211 & 433 & 343 & 565 & 157 & 379 & 121 & 506 & 228 & 475 & 67 & 289 & 199 & 416 & 13 & 360 & 577\\
 526 & 148 & 495 & 87 & 309 & 219 & 436 & 33 & 255 & 622 & 382 & 104 & 346 & 568 & 165 & 75 & 292 & 514 & 231 & 453 & 363 & 585 & 177 & 424 & 16\\
 95 & 312 & 534 & 126 & 498 & 258 & 605 & 222 & 444 & 36 & 571 & 168 & 390 & 107 & 329 & 239 & 456 & 53 & 300 & 517 & 402 & 24 & 366 & 588 & 185\\ \hline
 290 & 507 & 229 & 471 & 68 & 578 & 200 & 417 & 14 & 356 & 141 & 488 & 85 & 302 & 549 & 434 & 26 & 273 & 620 & 212 & 122 & 344 & 561 & 158 & 380\\
 454 & 71 & 293 & 515 & 232 & 17 & 364 & 581 & 178 & 425 & 310 & 527 & 149 & 491 & 88 & 623 & 220 & 437 & 34 & 251 & 161 & 383 & 105 & 347 & 569\\
 518 & 240 & 457 & 54 & 296 & 181 & 403 & 25 & 367 & 589 & 499 & 91 & 313 & 535 & 127 & 37 & 259 & 601 & 223 & 445 & 330 & 572 & 169 & 386 & 108\\
 57 & 279 & 521 & 243 & 465 & 375 & 592 & 189 & 406 & 3 & 538 & 135 & 477 & 99 & 316 & 201 & 448 & 45 & 262 & 609 & 394 & 111 & 333 & 555 & 172\\
 246 & 468 & 65 & 282 & 504 & 414 & 6 & 353 & 600 & 192 & 77 & 324 & 541 & 138 & 485 & 270 & 612 & 209 & 426 & 48 & 558 & 155 & 397 & 119 & 336\\ \hline
 441 & 38 & 260 & 602 & 224 & 109 & 326 & 573 & 170 & 387 & 297 & 519 & 236 & 458 & 55 & 590 & 182 & 404 & 21 & 368 & 128 & 500 & 92 & 314 & 531\\
 610 & 202 & 449 & 41 & 263 & 173 & 395 & 112 & 334 & 551 & 461 & 58 & 280  & 522 & 244 & 4 & 371 & 593 & 190 & 407 & 317 & 539 & 131 & 478 & 100\\
 49 & 266 & 613 & 210 & 427 & 337 & 559 & 151 & 398 & 120 & 505 & 247 & 469 & 61 & 283 & 193 & 415 & 7 & 354 & 596 & 481 & 78 & 325 & 542 & 139\\
  213 & 435 & 27 & 274 & 616 & 376 & 123 & 345 & 562 & 159 & 69 & 286 & 508 & 230 & 472 & 357 & 579 & 196 & 418 & 15 & 550 & 142 & 489 & 81 & 303\\
 252 & 624 & 216 & 438 & 35 & 570 & 162 & 384 & 101 & 348 & 233 & 455 & 72 & 294 & 511 & 421 & 18 & 365 & 582 & 179 & 89 & 306 & 528 & 150 & 492\\ \hline
 597 & 194 & 411 & 8 & 355 & 140 & 482 & 79 & 321 & 543 & 428 & 50 & 267 & 614 & 206 & 116 & 338 & 560 & 152 & 399 & 284 & 501 & 248 & 470 & 62\\
 11 & 358 & 580 & 197 & 419 & 304 & 546 & 143 & 490 & 82 & 617 & 214 & 431 & 28 & 275 & 160 & 377 & 124 & 341 & 563 & 473 & 70 & 287 & 509 & 226\\
 180 & 422 & 19 & 361 & 583 & 493 & 90 & 307 & 529 & 146 & 31 & 253 & 625 & 217 & 439 & 349 & 566 & 163 & 385 & 102 & 512 & 234 & 451 & 73 & 295\\
 369 & 586 & 183 & 405 & 22 & 532 & 129 & 496 & 93 & 315 & 225 & 442 & 39 & 256 & 603 & 388 & 110 & 327 & 574 & 166 & 51 & 298 & 520 & 237 & 459\\
 408 & 5  & 372 & 594 & 186 & 96 & 318 & 540 & 132 & 479 & 264 & 606 & 203 & 450 & 42 & 552 & 174 & 391 & 113 & 335 & 245 & 462 & 59 & 276 & 523\\ \hline
\end{array}
$$
}
\end{example}

\noindent More research into different properties and various examples can be found in the
thesis of the second author (\cite{thesis}). The reader is also referred to the website (\cite{website}).\\

\noindent {\bf Acknowledgment:} The authors like to thank Michiel de Bondt for various stimulating and helpful discussions.

\medskip\noindent
Harm Derksen\\
Dept. of Mathematics, Univ. of Michigan\\
East Hall, 530 Church Street\\
Ann Arbor, MI 48109-1043\\
\href{mailto:hderksen@umich.edu}{hderksen@umich.edu}\\

\noindent
Christian Eggermont\\
Dept. of Mathematics, Univ. of Nijmegen\\
6525 ED Nijmegen, The Netherlands\\
\href{mailto:C.Eggermont@science.ru.nl}{C.Eggermont@science.ru.nl}\\

\noindent
Arno van den Essen\\
Dept. of Mathematics, Univ. of Nijmegen\\ 
6525 ED Nijmegen, The Netherlands\\
\href{mailto:A.vandenEssen@science.ru.nl}{A.vandenEssen@science.ru.nl}\\

\end{document}